\documentclass[a4paper,10pt]{article}
\usepackage{theorem}
\usepackage{amsfonts}
\usepackage{amsmath} 
\usepackage{latexsym}
\usepackage{array}
\usepackage{amssymb}
\global
\theorembodyfont{\it}
\newtheorem{thm}{Theorem}[section]
\newtheorem{defi}[thm]{Definition}
\newtheorem{prop}[thm]{{\sc Proposition }}

\newtheorem{lem}[thm]{Lemma} 
\newtheorem{rem}[thm]{Remark}
\newtheorem{cor}[thm]{Corollary}
\newtheorem{acknowledgements}[thm]{Acknowledgements}
\theoremheaderfont{\sc}
\newtheorem{assumption}[thm]{Assumption}
\newcommand{\R}{  \mathbb{R}   }
\newcommand{\eps}{\varepsilon}
\newcommand{\e}{  \text{e}   }
\newcommand{\C}{  \mathbb{C}   }
\newcommand{\Z}{  \mathbb{Z}   }
\newcommand{\N}{  \mathbb{N}   }
\newcommand{\Q}{  \mathbb{Q}   }
\newcommand{\dem}{\noindent \textit{Proof: }}
\newcommand{\findem}{ }

\numberwithin{equation}{section} 
\newtheorem{notations}[thm]{Notations}

\date{}

\begin{document} 

\title{\vskip -1cm The WKB method and geometric instability for non linear Schr\"odinger equations on surfaces} 
 \author{ Laurent Thomann}

%\email{laurent.thomann@math.u-psud.fr}
\maketitle

\section{Introduction}

 \indent  Let $(M,g)$ be a  Riemannian surface (i.e. a Riemannian manifold of dimension $2$),
 orientable or not. We assume that $M$ is either compact or a compact perturbation of the euclidian space, so that the Sobolev embeddings are true. Consider $\Delta=\Delta_g$ the Laplace-Beltrami operator.
 In this paper we are interested in constructing WKB approximations for the non linear cubic Schr\"odinger equation
% $(NLS)$

\begin{equation}\label{nls}
\left\{
\begin{array}{l}
i \partial_t u(t,x)+\Delta u(t,x)   =  \eps|u|^{2} u(t,x), \quad \eps=\pm1\\
u(0,x)= u_0(x) \in H^{\sigma}(M)
\end{array}
\right.
\end{equation}
that is, given a small parameter $0<h<1$ and an integer $N$, functions $u_N(h)$ satisfying
\begin{equation}
i \partial_t u_N(h)+\Delta u_N(h)   =  \eps|u_N(h)|^{2} u_N(h)+R_N(h)
\end{equation}
with $\|u_N(h)\|_{H^{\sigma}} \sim 1$ and $\|R_N(h)\|_{H^{\sigma}}\leq C_N h^N.$\\
Here $h$ is introduced so that $u_N(h)$ oscillates with frequency $\sim\frac1h$.\\
These approximate solutions to $(\ref{nls})$ will lead to some instability properties in the following sense 
(where $h^{-1}$ will play the role of $n$):
\begin{defi}\label{definstab} We say that the Cauchy problem $(\ref{nls})$ is unstable near $0$ in
 $H^{\sigma}(M)$, if for all  $C>0$ there exist times  
$t_n \longrightarrow0$ and $u_{1,n}, u_{2,n} \in H^{\sigma}(M)$ solutions of (\ref{nls}) so that
\begin{eqnarray}
 \| u_{1,n}(0)\|_{H^{\sigma}(M)}, \;  \| u_{2,n}(0)\|_{H^{\sigma}(M)} &\leq& C   \nonumber \\
 \| u_{1,n}(0)- u_{2,n}(0)\|_{H^{\sigma}(M)} &\longrightarrow& 0 \nonumber\\ 
 \limsup \| u_{1,n}(t_n)- u_{2,n}(t_n)\|_{H^{\sigma}(M)} &\geq& \frac{1}{2}C  \nonumber
\end{eqnarray}
when $n \longrightarrow+\infty$.
\end{defi}
\noindent This means that the problem is not uniformly well posed, if we refer to the following definition:
\begin{defi}
Let $\sigma \in \R$. Denote by $B_{R,\sigma}$ the ball of radius $R$ in $H^{\sigma}$. We say that the Cauchy problem $(\ref{nls})$ if uniformly well posed in $H^{\sigma}$ if the flow map
$$u_0 \in B_{R,\sigma}\cap H^1(M) \longmapsto \Phi_t(u_0)\in H^{\sigma}(M)$$
is uniformly continuous for any $t$.
\end{defi}
We now state our instability result:
\begin{prop}\label{corinstab}
Let $0<\sigma<\frac{1}{4}$, and assume that $M$ has a stable and non degenerated
 periodic geodesic (see Assumptions 
\ref{assump1} and \ref{assump2} ), then the Cauchy problem $(\ref{nls})$ is not uniformly well-posed.
\end{prop}

\noindent This problem is motivated by the following results: In $\cite{BGT1}$, N. Burq, P. G\'erard and N. Tzvetkov show that 
$(\ref{nls})$ is unstable on the sphere $\mathbb{S}^2$ for $0<\sigma<\frac{1}{4}$. 
In fact they construct solutions of \eqref{nls} of the form
\begin{equation}\label{ans}
u_n^{\kappa}(t,x)=\kappa \e^{i\lambda_n^{\kappa}t}(n^{\frac14-\sigma}\psi_n(x)+r_n(t,x))
\end{equation}
where $0<\kappa<1$, $\psi_n=(x_1+ix_2)^n$ is a spherical harmonic which concentrates on the equator of the sphere 
when $n\longrightarrow+ \infty$ and where $r_n$ is an error term which is small. To obtain instability, they consider 
$\kappa_n\longrightarrow \kappa$, then 
\begin{equation*}
 \| u_{1,n}(0)- u_{2,n}(0)\|_{H^{\sigma}(\mathbb{S}^2)} \lesssim |\kappa'-\kappa|\longrightarrow 0 
\end{equation*}
but
\begin{equation*}
 \| u^{\kappa}_{n}(t_n)- u^{\kappa'}_{n}(t_n)\|_{H^{\sigma}(\mathbb{S}^2)} \gtrsim \kappa
|\e^{i\lambda_n^{\kappa}t_n}-\e^{i\lambda_n^{\kappa'}t_n}|\longrightarrow 2\kappa
\end{equation*}
with a suitable choice of $t_n\longrightarrow 0$. \\
We follow this strategy but as the surface is not rotation invariant, the ansatz will be 
more complicated than \eqref{ans}.\\
This result is sharp, because in $\cite{BGT2}$ they show that  $(\ref{nls})$ is uniformly well posed on  $\mathbb{S}^2$ when $\sigma>\frac{1}{4}$.\\
On the other hand, in $\cite{Bourgain}$ J. Bourgain shows that $(\ref{nls})$ is  uniformly well posed on the rational thorus $\mathbb{T}^2$ when $\sigma>0$.\\
%% therefore instability is strongly related to the geometry of $M$.
These results show how the geometry of $M$ can lead to instability fot the equation \eqref{nls}. Therefore it seems reasonable to obtain a result like Proposition \ref{corinstab} with purely geometric assumptions.
\\[10pt]
We first make the following assumption on $M$:

\begin{assumption}\label{assump1}
 The manifold $M$ has a periodic geodesic.
%$$\lambda \not \in \pi\mathbb{Q}$....Floquet, Poincar\'e...
\end{assumption}
Denote by $\gamma$ such a geodesic, then there exists a system of coordinates $(r,s)$ near $\gamma$, 
say for $(r,s)\in]-r_0, r_0[\times \mathbb{S}^1$, called Fermi coordinates such that (see \cite{Klingenberg}, p. $80$)
\begin{enumerate}
\item The curve $r=0$ is the geodesic $\gamma$ parametrized by arclength and 
\item The curves $s=$ constant are geodesics parametrized by arclength. 
The curves $r=$ constant meet these curves orthogonally .
\item In this system the metric writes
\[g=\left(\begin{array}{cc} 
1 & 0 \\ 0 & a^2(s,r) 
\end{array} \right). \]
\end{enumerate}
We set the length of $\gamma$ equal to $2\pi$. Denote by $R(r,s)$ the Gauss curvature at $(r,s)$, then 
$a$ is the unique solution of 
\begin{equation}\label{equationa}
\left\{
\begin{array}{l}
\frac{\partial^2a}{\partial r^2}+R(r,s)a  = 0\\[8pt]
a(0,s)=1,\frac{\partial a}{\partial r}(0,s)=0.
\end{array}
\right.
\end{equation}
The initial conditions traduce the fact that the curve $r=0$ is a unit-speed geodesic.
In these coordinates  the Laplace-Beltrami operator is
%\begin{eqnarray}
$$\Delta := \frac{1}{\sqrt{ \text{det}g} } \text{div}(
\sqrt{ \text{det}g}\; g^{-1}\nabla)=\frac{1}{a} \partial_s(\frac{1}{a}\partial_s)+\frac{1}{a} \partial_r(a\partial_r).$$
%\end{eqnarray}
A function on $M$, defined locally near $\gamma$, can be identified with a function of 
$[0,2\pi]\times ]-r_0,r_0[$ such that
$$\forall (s,r)\in [0,2\pi] \times ]-r_0,r_0[\quad f(s+2\pi,r)=f(s,\omega r)$$
where $\omega=1$ if $M$ is orientable and $\omega=-1$ if $M$ is not. Define
\begin{equation}\label{defomega_1}
\omega_1=\frac12(\omega-1)\in \{-1,0\}.
\end{equation}
From  $(\ref{equationa})$ we deduce that $a$ admits the Taylor expansion 
\begin{equation}\label{dvpta}
a=1-\frac12R(s)r^2+R_3(s)r^3+\cdots+R_p(s)r^p+o(r^p),
\end{equation}
and as $a(s+2\pi,r)=a(s,\omega r)$, we deduce $R(s+2\pi)=R(s)$ and for all $j\geq 3$,  $R_j(s+2\pi)=\omega^j R_j(s)$.

\noindent Let $p_2=\frac{1}{a^2}\sigma^2+\rho^2$ be the principal symbol of $ \Delta $, and 
\begin{equation}\label{hamiltonien}
\left\{
\begin{array}{ll}
\frac{\text{d}}{\text{d}t}{s}(t)=\frac{\partial p_2}{\partial \sigma}=\frac{2\sigma}{a^2},&\frac{\text{d}}{\text{d}t}{\sigma}(t)=-\frac{\partial  p_2}{\partial s}
=-\partial_s({\frac{1}{a^2}})\sigma^2\\[3pt]
\frac{\text{d}}{\text{d}t}{r}(t)=\frac{\partial  p_2}{\partial \rho}=2\rho,&\frac{\text{d}}{\text{d}t}{\rho}(t)=-\frac{\partial  p_2}{\partial r}=-\partial_r{(\frac{1}{a^2}})\sigma^2\\[3pt]
s(0)=s_0,\;\sigma(0)=\sigma_0,&r(0)=r_0,\;\rho(0)=\rho_0,
\end{array}
\right.
\end{equation}
its associated hamiltonian system, where $p_2=p_2(s(t),r(t),\sigma(t),\rho(t))$. 
The system \eqref{hamiltonien} admits a unique solution and defines the hamiltonian fow
$$\Phi_t:(s_0,\sigma_0,r_0,\rho_0)\longmapsto (s(t),\sigma(t),r(t),\rho(t)).$$
The curve 
$\Gamma=\{(s(t)=t,\sigma(t)=1/2,r(t)=0,\rho(t)=0), t\in [0, 2\pi]\}$ is solution of $(\ref{hamiltonien})$ and its projection in the $(r,s)$ space 
is the curve $\gamma$. Now denote by $\phi$ the Poincar\'e map associated to 
the trajectory $\Gamma$ and to the hyperplane $\Sigma=\{s=0\}$. There exists a neighborhood $\mathcal{N}$ of 
$(\sigma=1/2,r=0,\rho=0)$ such that the following makes sense:
solve the system $(\ref{hamiltonien})$ with the initial conditions $(0,\sigma_0,r_0,\rho_0)\in \{0\}\times\mathcal{N}$ 
and let $T$ be such that $s(T)=2\pi$, then $\phi$ is the application 
$$\phi:(r_0, \rho_0)\longmapsto  (r(T), \rho(T)). $$
Moreover, the Poincar\'e map is continuously differentiable (see $\cite{Perko}$ p. $193$).
To obtain its differential  $\text{d}\phi(0,0)$ at $(0,0)$, we linearize the system  
 $(\ref{hamiltonien})$ about the orbit $\Gamma$, i.e.
\begin{equation}\label{linhamilton}
\left\{
\begin{array}{ll}
\frac{\text{d}}{\text{d}t}{s}(t)=2\sigma,&\frac{\text{d}}{\text{d}t}{\sigma}(t)=0\\[2pt]
\frac{\text{d}}{\text{d}t}{r}(t)= 2\rho,&\frac{\text{d}}{\text{d}t}{\rho}(t)=-\frac12 R(s(t)) r,
\end{array}
\right.
\end{equation}
then $\sigma=\frac{1}{2}$, $s(t)=t$ and 
\begin{equation}\label{systdiff}
%\begin{displaymath}
\frac{\text{d}}{\text{d}t}{\left(\begin{array}{c} 
{r} \\ {\rho} 
\end{array} \right)}=\left(\begin{array}{cc} 
0 & 2 \\ - R/2 & 0 
\end{array} \right){\left(\begin{array}{c} 
{r} \\ {\rho} 
\end{array} \right)} .
%\end{displaymath}
\end{equation}
Hence the application $\text{d}\phi(0,0)$ is 
$$\text{d}\phi(0,0):(r_0, \rho_0)\longmapsto  (r(2\pi), \rho(2\pi)), $$
where $(r,\rho)$ solve \eqref{systdiff}. As $\text{d}\phi(0,0)$ is symplectic, it admits 
two eigenvalues $\mu$ and $\mu^{-1}$ that are called the characteristic multipliers of the 
system \eqref{systdiff}.
We add the following assumption on $\gamma$, which can be formulated in terms of the eigenvalues 
of $\text{d}\phi(0,0)$:
\begin{assumption}\label{assump2}
The geodesic $\gamma$ is stable, i.e. $\text{d}\phi(0,0)$ is a rotation. Then 
the multipliers take the form $\mu=\e^{i\lambda}$ and $\mu^{-1}=\e^{-i\lambda}$ 
with  $\lambda\in  \mathbb{R}$. We assume moreover that 
%$\lambda \not \in \pi \mathbb{Q}$.
there exist $\tau,\mu>0$ such that
\begin{equation}\label{hyp}
\forall (p,q)\in \Z\times \N\quad |p-q\frac{\lambda}{\pi}|\geq \frac{\mu}{|(p,q)|^{\tau}},
\end{equation}
where $|(p,q)|=|p|+|q|$. When this condition is fulfilled, we say that $\gamma$ is non degenerated.
\end{assumption}

\begin{rem} Almost every $\lambda\in \R$ satisfies \eqref{hyp} with $\tau>1$. This is an easy 
consequence of \cite{AlinhacGerard} p. $159$, e.g.
\end{rem}
\noindent Notice that the function $r$ which satisfies \eqref{systdiff} is solution of 
\begin{equation}\label{Eqa_0}
\ddot{y}(s)+R(s)y(s)=0.
\end{equation}
Consider $a_0$ the solution of \ref{Eqa_0} with initial conditions $a_0(0)=1$ and $\dot{a}_0(0)=i$. Then, from 
the Floquet theory, there exists a $2\pi$-periodic function $P$ so that 
$$a_0(s)=\e^{i\frac{\lambda}{2\pi}s}P(s)$$
(or  $a_0(s)=\exp{(-i\frac{\lambda}{2\pi}s})P(s)$, but $\lambda$ can be replaced with $-\lambda$).\\
Here, and in all the paper we denote by $\dot{f}=\frac{\text{d}}{\text{d}s}f$ if $f$ is differentiable. 
This notation is motived by the fact that $s$ will play the role of a time variable (see section $\ref{sectionWKB}$).
\\[10pt]
\noindent In order to prove Proposition \ref{corinstab}, we construct stationnary approximate solutions 
of \eqref{nls}, as stated in the following Theorem
\begin{thm}\label{thmWKB}
Assume $\ref{assump1}$ and $\ref{assump2}$. Let $h\in ]0,1]$ such that $\frac1h\in\N$, 
let $\kappa, \sigma>0$ and $k\in \N$.
Let $\lambda$ be given by Assumption \ref{assump2} and  $\omega_1$ by \eqref{defomega_1}. \\
Define $E_0(k)=-\frac{1}{4\pi}\lambda 
+\frac12k(\omega_1-\frac{\lambda}{\pi})$ and $\delta =\kappa h^{\sigma}$.\\
Then for all $N\in\N$, there exist $\lambda_{N}(k)\in \R$ and a family $u_N(h)$ such that \\
$C_1\delta \leq\|u_N(h)\|_{L^2(M)}\leq C_2\delta$ with $C_1,C_2>0$ independent of $N$ and $h$, and
\begin{equation}\label{eqthm}
-\Delta u_N(h) = \lambda_N(k) u_N(h) -\eps |u_N(h)|^2u_N(h)+h^Ng_N(h) 
\end{equation}
with for all $N \in \N$
$$\|h^Ng_N(h)\|_{H^n(M)}\lesssim h^{N-n}.$$
Moreover $$\lambda_N(k)=\frac{1}{h^2}-\frac{2}{h}E_0(k)+
\frac{1}{\sqrt{h}}\eps\delta^2C_0+\mathcal{O}(1),$$
where $C_0>0$ is independent of $\eps$ and $\delta$.
\end{thm}

\begin{rem}
The analog of Theorem $\ref{thmWKB}$ was  proved by J. Ralston in $\cite{Ralston1}$ for the linear 
case $(\eps=0)$, with the same type of 
assumptions.
\end{rem}

\begin{rem}
Consider the more general equations 
\begin{equation}\label{nlsF}
i \partial_t u+\Delta u   = F(u),
\end{equation}
where $F: \C \longrightarrow \C $ is a $\mathcal{C}^{\infty}$ function.
The result of  Theorem $\ref{thmWKB}$ is likely to hold with other nonlinearities $F(u)$, for 
example for $F(z)=z^3$, $F(z)=z^4$ or $F(z)=(1+|z|^2)^{\alpha}z$ with $\alpha<1$. However, 
the instability phenomenon is strongly related to the gauge invariance of the equation \eqref{nlsF}.
% For instance 
%X proves in \cite{X} that, in the case  $F(z)=z^3$,  \eqref{nlsF} is well-posed in $H^{-\frac12}(\mathbb{S}^2)$; 
%but we expect Corollary \ref{corinstab} to hold in the case $F(z)=(1+|z|^2)^{\alpha}z$.
\end{rem}

\begin{rem}
The restriction \eqref{hyp} seems to be  purely technical. 
For instance, when $M=\mathbb{S}^2$, $\lambda=0$ but the 
result holds (see \cite{BGT1}).
\end{rem}

\noindent The scheme of the paper is the following: Thanks to a scalling, we reduce the problem $(\ref{eqthm})$ to the resolution of linear Schr\"odinger equations with a harmonic time dependent potential, and we will see, using  Assumption  $\ref{assump2}$, that these equations have periodic solutions.
 To prove Proposition $\ref{corinstab}$ we show that the family $u_N(h)$ provide good approximations of $(\ref{nls})$ in times where instability occurs.

%We work in Fermi coordinates $(s,r)$, with $s\in [0,2\pi]$ and $r\in ]-\epsilon,\epsilon[$.
%\begin{equation}
%\Delta=\frac{1}{a} \partial_s(\frac{1}{a}\partial_s)+\frac{1}{a} \partial_r(a\partial_r),
%\end{equation}
%where $a$ is given by
%The function $R(s,r)>0$ is the curvature of $M$ at the  point  $(s,r)$ and is 
%$2\pi$-periodic in $s$.
%$$R(s,r)=R(s,0)+\frac{1}{2}\frac{\partial^2R}{\partial r^2}(s,0)r^2+o(r^2),$$
%avec $\frac{1}{2}\frac{\partial^2R}{\partial r^2}(s,0)<0$.
%$$\frac{\partial^2a}{\partial r^2}(s,0)=-R(s,r)<0, \, \frac{\partial^3a}{\partial r^3}(s,0)=0 \,
%\text{et} \,\frac{\partial^4a}{\partial r^4}(s,0)>0.  $$
%For this metric, the Laplace-Belrami operator writes
%$$\Delta=\frac{1}{a} \partial_s(\frac{1}{a}\partial_s)+\frac{1}{a} \partial_r(a\partial_r).$$

\begin{notations}
In this paper $c$, $C$ denote constants the value of which may change from line to line. We use the notations $a\sim b$,  
$a\lesssim b$ if $\frac1C b\leq a\leq Cb$, $a\leq Cb$ respectively. By $\delta_{i,j}$ we mean the Kronecker symbol, i.e. 
$\delta_{i,j}=0$ for $i\neq j$ and $\delta_{i,i}=1$.
\end{notations}

\begin{rem}
In the sequel we do not always mention the dependence on $h$ of the functions: we will write $u$, $f$, $r_i$, $\dots$ 
instead of $u_h$, $f_h$, $r_{i,h}$, $\dots$
\end{rem}
\begin{acknowledgements}
The author would like to think his adviser N. Burq for his constant guidance in this work, P. Pansu for his help 
in the frame of geometry, and B. Helffer for having pointed out the reference \cite{Combescure}.
\end{acknowledgements}

\section{The WKB construction}\label{sectionWKB}
Consider the equation 
\begin{equation}\label{elliptique}
-\Delta u = \lambda u -\eps |u|^2u.
\end{equation}
Given $h>0$, we are looking for a solution of the form 
\begin{equation}\label{ansatzu}
u=\delta h^{-\frac{1}{4}}\e^{i\frac{s}{h}}f(s,r,h)
\end{equation}
 where $\delta=\kappa h^{\sigma}$, with $\kappa>0$ and $0\leq \sigma\leq \frac14$.
We try to find a solution of \eqref{elliptique} of the form $u\sim \sum_{j\geq 0}h^{j/2}u_j$. Thus 
we write $\lambda\sim h^{-2}\sum_{j\geq 0}h^{j/2}\lambda_j$. As we will see, identifying each power of $h$ will lead to a linear equation which can be solved with a suitable choice of $\lambda_j$.\\
 Choose $h$ such that $h^{-1}\in \N$, this ensures 
that $\exp{i\frac{s}{h}}$ is $2\pi-$periodic. Such a condition on $h$ is natural and is known as a Bohr-Sommerfeld quantification condition.\\
 With the ansatz \eqref{ansatzu}, equation  $(\ref{elliptique})$ becomes 
\begin{eqnarray}
-\frac{1}{a^2}(\frac{2i}{h}\partial_{s}f+\partial^2_sf-\frac{1}{h^2}f) -\frac1a\partial_{s}
(\frac1a)(\frac ihf+\partial_sf) \nonumber \\
-\partial_r^2f-\frac{\partial_ra}{a}\partial_rf=\lambda f -\eps 
\delta^2h^{-\frac{1}{2}}|f|^2f. \label{equationf}
\end{eqnarray}
We make the change of variables $x=\frac{r}{\sqrt{h}}$ and set $v(s,x,h)=f(s,\sqrt{h}x,h)$. Thus 
$\partial_rf=\frac{1}{\sqrt{h}}\partial_xv$ and 
$\partial^2_rf=\frac{1}{{h}}\partial^2_xv$. Therefore we now have to find
 $v\sim \sum_{j\geq 0}h^{j/2}v_j$. Using $(\ref{dvpta})$ we obtain the following Taylor expansions in $h$
%\begin{eqnarray*}
$$\frac{1}{a^2}=1+hRx^2-2h^{\frac32}R_3x^3+\mathcal{O}(h^2),$$
$$a^{-1}\partial_{s}(a^{-1})=\mathcal{O}(h)\quad \text{and}\quad  
a^{-1}\partial_ra =\mathcal{O}(h^{\frac12}).$$ 
%\end{eqnarray*}
Equation $(\ref{equationf})$ can therefore be written, after multiplication by $\frac12h$
\begin{eqnarray}\label{equationv}
i\partial_{s}v+\frac12\partial^2_xv    -\frac12Rx^2v= 
\frac{1-\lambda h^2}{2h}v+h^{\frac12}R_3x^3v   +\frac12\eps \delta ^2 h^{\frac12} |v|^2v 
+hPv
\end{eqnarray}
where 
\begin{equation}\label{defP}
P=A_1\partial^2_{s}+ A_2\partial_{s}+ A_3\partial_{r}+A_4
\end{equation}
is a second order differential operator with coefficients $A_j=A_j(s,x,h)$ satisfying 
$A_j(s+2\pi,x,h)=A_j(s,\omega x,h)$ for  $0\leq j \leq 4$.\\
Denote by $E=\frac{1-\lambda h^2}{2h}=E_0+h^{\frac12}E_1+\cdots+h^{\frac{p}2}E_p+o(h^{\frac{p}2})$ and 
write $v=v_0+h^{\frac12}v_1+\cdots+h^{\frac{p}2}v_p+o(h^{\frac{p}2})$ and 
by identifying the powers of $h$ we obtain the system of equations:
%\newcounter{\toto}
%\newcounter{toto} \setcounter{\theequation}{toto.a}
\begin{eqnarray}\label{system}
\left(i\partial_{s}+\frac12\partial^2_x    -\frac12Rx^2-E_0\right)v_0 &=& 0 \label{Eqv_0}\\
\left(i\partial_{s}+\frac12\partial^2_x    -\frac12Rx^2-E_0\right)v_1 &=&    E_1 v_0+R_3x^3v_0
 +\frac12\eps \delta ^2 |v_0|^2v_0 \label{Eqv_1}\\
\cdots &=& \cdots  \nonumber \\
\left(i\partial_{s}+\frac12\partial^2_x    -\frac12Rx^2-E_0\right)v_p &=&
+E_p v_0+Q_{p}. \label{Eqv_p}\\
\cdots &=& \cdots  \nonumber 
\end{eqnarray}
so that the $(j+1)$th equation of unknown $(v_j,E_j)$ correspond to the annihilation of the coefficient 
of $h^{\frac{j}{2}}$ in $(\ref{equationv})$.\\
Here $Q_{p}$ is a function which only depends on $x,s$, $(v_j)_{j\leq p-1}$ and $(E_j)_{j\leq p-1}$.

\begin{rem}
Notice that thanks to the scalling, we have reduced the problem $(\ref{elliptique})$ to the resolution of  
linear equations. However we have to solve them exactly; no smallness assumption on $x$ is possible, as $x$ 
can be of size $\sim \frac{1}{\sqrt{h}}$.
\end{rem}
Once we have solved the previous equations, we will be able to construct approximate solutions of 
\eqref{elliptique}; more precisely, we will obtain the following proposition, which is the main result of this section.

\begin{prop}\label{propWKB}
Let $\chi \in \mathcal{C}_0^{\infty}(]-r_0, r_0[)$ be such that $0\leq \chi \leq1$, $\chi=1$ on $[-r_0/2, r_0/2]$ 
and suppose moreover that $\chi$ is an even function.  
Denote by 
\begin{equation}\label{expru_p}
u_p=\delta h^{-\frac14} \chi(r)\e^{i\frac{s}{h}}(v_0+h^{\frac12}v_1+\cdots
+h^{\frac{p}{2}}v_p)(s,{\frac{r}{\sqrt{h}}})
\end{equation}
and by
\begin{equation}\label{exprlambda_p}
\lambda_p=\frac{1}{h^2}-\frac{2}{h}(E_0+h^{\frac12}E_1+\cdots+h^{\frac{p}{2}}E_p).
\end{equation}
Then $u_p$ satisfies $\|u_p\|_{L^2(M)}\sim \delta$ and 
\begin{equation}\label{exprelliptique_p}
-\Delta u_p = \lambda_p u_p -\eps |u_p|^2u_p+ h^{\frac{p-1}{2}}g_p(h) 
\end{equation}
with 
$$\forall h\in ]0,1],\;\forall n\in \N,\quad \| h^{\frac{p-1}{2}}g_p(h) \|_{H^n} \lesssim \delta h^{\frac{p-1}{2}-n}.$$
\end{prop}

\subsection{Preliminaries: the analysis of the linear equations}
We will solve the system $(\ref{system})$ for $x\in \R$, even if  the function $a$ is only defined 
for $|r|\leq r_0$ i.e. for $x \leq \frac{r_0}{\sqrt{h}}$, but the equations $(\ref{system})$ make sense for 
$x\in \R$. At the end we will use a cutoff argument to construct a proper function on $M$.

\noindent Consider the Hilbertian basis of $L^2(\R)$ composed of the Hermite functions $(\varphi_k)_{k\geq 0}$ 
which are the eigenfunctions of the harmonic oscillator $H=-\frac12\partial^2_x+\frac12x^2$, i.e. 
$H\varphi_k=(k+\frac12)\varphi_k$. Moreover $\varphi_k(x)=P_k(x)\e^{-x^2/2}$ where $P_k$ is a polynomial of degree $k$ 
with $P_k(-x)=(-1)^kP_k(x)$. The link between the $s$-dependent operator $-\frac12\partial^2_x+\frac12R(s)x^2$ and $H$ 
is given by the following result proved by M. Combescure in $\cite{Combescure}$.

\begin{thm}\label{thmcombescure}
Let $a_0:\R\longrightarrow \C$ be the solution of $(\ref{Eqa_0})$ with $a_0(0)=1$, $\dot{a}_0(0)=i$. Define 
$$\alpha=\log{|a_0|},\;\beta=\frac{1}{2i}\log{\frac{a_0}{\overline{a_0}}},$$
let the unitary transform $T(s)$ be defined by 
$$T(s)=\e^{i\dot{\alpha}(s)x^2/2}\e^{-i\alpha(s)D},\;\text{where}\;D=-\frac{i}{2}(x\cdot\nabla+\nabla\cdot x),$$
and let $U(s,\tau)$ be the unitary evolution operator for
$-\frac12\partial^2_x+\frac12R(s)x^2$, i. e. $U(s,\tau)\varphi$ is the unique solution of the problem 
\begin{equation*}
\left\{
\begin{array}{l}
(i \partial_s +\frac12\partial^2_x-\frac12R(s)x^2)u=0\\
u(\tau,x)= \varphi(x) \in L^2({\R})
\end{array}
\right.
\end{equation*}
Then we have for any $s, \tau\in \R$
$$U(s,\tau)=T(s)\e^{-i(\beta(s)-(\beta(\tau))H}T(\tau)^{-1}.$$
\end{thm}

\begin{rem}
The functions $\alpha$ and $\beta$ are well defined: suppose that there exists $s_0$ such 
that $a_0(s_0)=0$, then $\text{Re} (a_0)$ and $\text{Im} (a_0)$ are linearly dependent, 
which is impossible with this choice of the initial conditions. 

\end{rem}

\begin{rem} Define $\theta(s)=\beta(s)-\frac{\lambda}{2\pi}s$ where $\lambda$ is given 
by Assumption \ref{assump2}. Then $\alpha$ and $\theta$ are $2\pi$-periodic real functions.
Moreover $\alpha(0)=\dot{\alpha}(0)= \beta(0)=\theta(0)=0$.
\end{rem}
Denote by $\mathcal{S}(\R)$ the Schwartz space, i.e. the space of smooth functions which are fast decreasing and 
their derivatives too.

\begin{prop}\label{propregularite}
Let $\psi_0\in \mathcal{S}(\R)$ and $E\in \C$. Let 
$f\in \mathcal{C}^{\infty}\left([0,2\pi]\times \R,\R\right)$ be such that 
$$\forall n \in \N, \;\forall s \in [0,2\pi],\quad \partial_s^nf(s,\cdot)\in \mathcal{S}(\R),$$
in other words $f \in \mathcal{C}^{\infty}\left([0,2\pi],\mathcal{S}(\R)\right)$.\\
Let $\psi \in  \mathcal{C}^1\left([0,2\pi],L^2(\R)\right)\cap \mathcal{C}^0\left([0,2\pi],H^2(\R)\right)$ 
be the solution of 

\begin{equation}\label{Eqpsi}
\left\{
\begin{array}{l}
i \partial_s\psi +\frac12\partial^2_x\psi -\frac12Rx^2\psi-E\psi=f\\
\psi(0,x)= \psi_0(x).
\end{array}
\right.
\end{equation}
Then  $\psi\in  \mathcal{C}^{\infty}\left([0,2\pi],\mathcal{S}(\R)\right)$.
% and  $$\forall n \in \N, \;\forall s \in \R,\quad \partial_s^nf(s,\cdot)\in \mathcal{S}(\R).$$
\end{prop}

\dem
By replacing $\psi$ with $\e^{iEt}\psi$, we can assume that $E=0$. The solution of equation $(\ref{Eqpsi})$ is given by 
\begin{eqnarray}\label{duhamel}
\psi(s,\cdot)&=&U(s,0)\psi_0-i\int_{0}^sU(s,\tau)f(\tau,\cdot)\text{d}\tau\nonumber\\
&=&T(s)\e^{-i\beta(s)H}\left(\psi_0 -i\int_{0}^se^{i\beta(\tau)H}T(\tau)^{-1}f(\tau,\cdot)\text{d}\tau\right).
\end{eqnarray}
As $D$ is a transport operator, we have
$$T,T^{-1} :\mathcal{C}^{\infty}\left([0,2\pi],\mathcal{S}(\R)\right) 
\longrightarrow \mathcal{C}^{\infty}\left([0,2\pi],\mathcal{S}(\R)\right),  $$
we only have to show that 
$$\e^{i\beta H} :\mathcal{C}^{\infty}\left([0,2\pi],\mathcal{S}(\R)\right) 
\longrightarrow \mathcal{C}^{\infty}\left([0,2\pi],\mathcal{S}(\R)\right). $$
This follows from the fact that  $\beta$ is regular and $\e^{iH}:\mathcal{S}(\R) \longrightarrow \mathcal{S}(\R)$.
\findem
\\[10pt]
The description of $U$ given in Theorem $\ref{thmcombescure}$ yields the following representation of $U(s,0)\varphi_k$:

\begin{prop}\label{propcalcul}
For all $k\in \N$ and  $s,x\in \R$ we have
\begin{equation}
U(s,0)\varphi_k(x)=\e^{i\dot{\alpha}(s)x^2/2}\e^{-i(\frac12+k)\beta(s)}\e^{-\frac12\alpha(s)}\varphi_k\left(x\e^{-\alpha(s)}\right).
\end{equation}
\end{prop}

\dem
According to Theorem $\ref{thmcombescure}$, and as $H\varphi_k=(k+\frac12)\varphi_k$,
$$U(s,0) \varphi_k = \e^{i\dot{\alpha}(s)x^2/2}\e^{-i(k+\frac12)\beta(s)}\e^{-i\alpha(s)D} \varphi_k. $$
Denote by $f(s)=\e^{-i\alpha(s)D} \varphi_k$. Then $f$ is solution of the transport equation 
\begin{equation*}
\partial_sf=-\frac12\dot{\alpha}(s)\big(x\partial_xf+\partial_x(xf)\big)
           =-\frac12\dot{\alpha}(s)\left(f+2x\partial_xf\right)
\end{equation*}
with Cauchy data $f(0,x)=\varphi_k(x)$.
Make the change of variables $\sigma=\alpha(s)$ and set $g(\sigma)=f(s)$. Therefore 
$g$ satisfies $\partial_{\sigma}g=-\frac12(g+2x\partial_xg)$. The equation 
$x=\dot{x}$, $x(0)=x_0$ admits the solution $x(\tau)=x_0\e^{\tau}$ and  the 
characteristics method gives $g(\tau,x(\tau))=\e^{-\frac12\tau}\varphi_k(x_0)
= \e^{-\frac12\tau}\varphi_k(x(\tau)\e^{-\tau})$, hence 
\begin{equation*}
f(s)=\e^{-\frac12\alpha(s)} \varphi_k(x\e^{-\alpha(s)}).
\end{equation*}
\findem

\begin{cor}\label{lemw_0}
Let $k \in  \N$, define $\omega_1=\frac12(\omega-1)$ and \\
$\mathcal{E}_{k_0}=E_0(k)=-\frac{1}{4\pi}\lambda 
+\frac12k(\omega_1-\frac{\lambda}{\pi})$. Then 
\begin{eqnarray}
w_k&=&\e^{-i s \mathcal{E}_{k_0}} U(s,0) \varphi_k \nonumber \\
   &=&\e^{-i s \mathcal{E}_{k_0}}\e^{i\dot{\alpha}(s)x^2/2}\e^{-i(\frac12+k)\beta(s)}\e^{-\frac12\alpha(s)}\varphi_k\left(x\e^{-\alpha(s)}\right) \label{expressionw_k}
\end{eqnarray}
is solution of the equation $(\ref{Eqv_0})$. 
\end{cor}

\dem
On the one hand, from Proposition $\ref{propcalcul}$ we deduce
\begin{eqnarray*}
w_k(s+2\pi,x)&=&\e^{-2i\pi\mathcal{E}_{k_0} }\e^{-i\lambda (\frac12+k)}w_k(s,x)=
\e^{-i k \omega_1 \pi}w_k(s,x)\\
&=& (-1)^{k \omega_1}w_k(s,x)=w_k(s,\omega x).
\end{eqnarray*}
On the other hand, $w_k$ satisfies $(\ref{Eqv_0})$ because of the definition of $U(s,0)$. 
\findem
\\[10pt]
Fix $k_0 \in \N$ and take $v_0=w_{k_0}$ with the previous choice of $\mathcal{E}_{k_0}=E_0(k_0)$. 
This choice corresponds to the $k_0$th level of energy for the harmonic oscillator.

\begin{rem}
Until now we didn't use the restriction \eqref{hyp}, but it will crucial in the following.
\end{rem}

\begin{prop}\label{propE_p}
For all $p \geq 0$, there exist $E_p \in \C$ and\\
 $v_p \in \mathcal{C}^{\infty}\left([0,2\pi],\mathcal{S}(\R)\right)$ 
which solve $(\ref{Eqv_p})$.
\end{prop}

\begin{rem}
As stated in Theorem $\ref{thmWKB}$, the $E_j$'s are in fact real numbers. 
This will be proved in Lemma $\ref{lemreel}$.
\end{rem}

\dem
We proceed by induction on $p\in \N$.\\
For $p=0$ the result was  proved in Corollary \ref{lemw_0}.\\
Let $p \geq 1$, and suppose that for all $j\leq p-1$ there exist $E_j\in \C$ and 
$v_j \in \mathcal{C}^{\infty}\left([0,2\pi],\mathcal{S}(\R)\right)$ 
which solve the $(j+1)$th equation of $(\ref{system})$. 
When $p\geq 2$, set  
\begin{eqnarray*}
\tilde{v}_{p-1} & = & h^{\frac12}v_1+\cdots +h^{\frac{p-1}{2}}v_{p-1},\\
\tilde{E}_{p-1} & = & h^{\frac12}E_1+\cdots +h^{\frac{p-1}{2}}E_{p-1}
\end{eqnarray*}
and  $\tilde{v}_{0}=\tilde{E}_{0}=0$.  The function $Q_{p}$ given by \eqref{Eqv_p} is the 
coefficient of $h^{\frac{p}{2}}$ in the expansion in h of
$$\tilde{E}_{p-1}\tilde{v}_{p-1} +\frac12\eps \delta ^2 |v_0+\tilde{v}_{p-1}|^2(v_0
+\tilde{v}_{p-1})+hP(v_0+\tilde{v}_{p-1}).$$ 
Now using the regularity of the $v_j$'s and the fact that $P$ defined by $(\ref{defP})$ is 
an operator
$$P:\mathcal{C}^{\infty}\left([0,2\pi],\mathcal{S}(\R)\right) 
\longrightarrow \mathcal{C}^{\infty}\left([0,2\pi],\mathcal{S}(\R)\right),$$
we obtain $Q_p \in \mathcal{C}^{\infty}\left([0,2\pi],\mathcal{S}(\R)\right). $\\
Moreover $Q_p$ satisfies, $\forall (s,x)\in [0,2\pi]\times \R$
$$Q_p(s+2\pi,x)=Q_p(s,\omega x)$$
because this property holds for the $v_j$'s,  $a$ and for the coefficients of $P$.\\
Define $F_p(s,x)=\e^{-i\dot{\alpha}(s)\e^{2\alpha(s)}x^2/2 }Q_p(s,xe^{\alpha(s)})$, then 
 $F_p \in \mathcal{C}^{\infty}\left([0,2\pi],\mathcal{S}(\R)\right)$ and satisfies 
$Q_p(s,x)=\e^{i\dot{\alpha}(s)x^2/2 }F_p(s,x\e^{-\alpha(s)})$ and 
$F_p(s+2\pi,x)=F_p(s,\omega x)$. Let us  decompose $F_p$ on the basis  
$(\varphi_j)_{ j\geq0}$: there exists a unique family of smooth functions 
$(g^p_j(s))_{ j\geq0} \in l^2(\N)$ so that 
\begin{equation}\label{defF}
F_p(s,y)=\sum_{j\geq0}g^p_j(s)\varphi_j(y).
\end{equation}
Then 
\begin{equation}\label{dvptQ_p}
Q_p(s,x)=\sum_{j\geq0}g^p_j(s)\e^{i\dot{\alpha}(s)x^2/2 }\varphi_j(x\e^{-\alpha(s)})
       = \sum_{j\geq0}h^p_j(s)w_j(s,x)
\end{equation}
where according to $(\ref{expressionw_k})$
\begin{equation}\label{defhj}
h^p_j(s)= \e^{i s E_0(j)}\e^{i(\frac12+j)\beta(s)}\e^{\frac12\alpha(s)} g^p_j(s).
\end{equation}
We have 
$$Q_p(s,\omega x)=  \sum_{j\geq0}h^p_j(s)w_j(s,\omega x),$$
but also
\begin{eqnarray*}
Q_p(s,\omega x)= Q_p(s+2 \pi,x)&=&  \sum_{j\geq0}h^p_j(s+2 \pi)w_j(s+2 \pi,x)\\
                       &=& \sum_{j\geq0}h^p_j(s+2 \pi)w_j(s,\omega x),
\end{eqnarray*}
and from the uniqueness of the $h^p_j$'s we deduce $h^p_j(s+2 \pi)=h^p_j(s)$.
%, for all $s \in [0,2\pi]$.

\noindent We are now looking for a solution of $(\ref{Eqv_p})$ of the form 
\begin{equation}\label{defvp}
v_p(s,x)=\sum_{j \geq 0} e^p_j(s)w_j(s,x)
\end{equation}
where the $e^p_j$'s are $2\pi$-periodic functions. 
%Indeed such a function satisfies
%\begin{equation*} 
%v_1(s+2\pi,x)=\sum_{j \geq 0} e_j(s+2\pi)w_j(s+2\pi,x)
%=\sum_{j \geq 0} e_j(s)w_j(s,\omega x)=v_1(s,\omega x).
%\end{equation*}
For all $j\geq 0 $, by Lemma \ref{lemw_0} we have 
$$\left( i\partial_{s}+\frac12\partial^2_x    -\frac12Rx^2 \right)\Big(e^p_j w_j\Big)=
i \dot{e}^p_jw_j+\big( \mathcal{E}_{k_0}- \mathcal{E}_{j} \big)e^p_j w_j,$$ 
hence we have to solve the equations 
\begin{equation}\label{Eqe_j}
i \dot{e}^p_j+\big( \mathcal{E}_{k_0}- \mathcal{E}_{j} \big)e^p_j=h^p_j+\delta_{j,k_0}E_p.
\end{equation}
As $ \mathcal{E}_{k_0}- \mathcal{E}_{j}= E_0(k_0)- E_0(j)=\frac12(k_0-j)(\omega_1-\frac{\lambda}{\pi})$, the solutions 
of \eqref{Eqe_j} take the form
\begin{equation}\label{formulee_j}
 e^p_j(s)=\e^{\frac12i(k_0-j)(\omega_1-\frac{\lambda}{\pi})s} 
\left(C^p_j-i\int_{0}^{s}h^p_j(\tau)\e^{-\frac12i(k_0-j)
(\omega_1-\frac{\lambda}{\pi})\tau}\text{d}\tau\right) 
\end{equation}
 for $j\neq k_0$, and 
$$e^p_{k_0}(s)=C^p_{k_0}-i\int_{0}^{s}h^p_{k_0}(\tau)\text{d}\tau -i E_ps.$$
The constants  $C^p_j\in \C $ and $E_p \in \C$ have  to be determined such that $e^p_j(s+2\pi)=e^p_j(s)$.\\
$\bullet$ Case $j=k_0$: 
$$e^p_{k_0}(s+2\pi)=-i\int_{0}^{2\pi}h^p_{k_0}(\tau)\text{d}\tau -2\pi i E_p +e^p_{k_0}(s),$$
thus $e^p_{k_0}$ is $2\pi$-periodic iff 
\begin{equation}\label{formuleE_p}
E_p=-\frac{1}{2\pi}\int_{0}^{2\pi}h^p_{k_0}\text{d}\tau.
\end{equation}
$\bullet$ Case $j\not =k_0$: \\
Denote by $\tilde{h}^p_j:\tau \longmapsto h^p_j(\tau) 
\e^{-i\frac12(k_0-j)(\omega_1-\frac{\lambda}{\pi})\tau}$ and by 
$K=\e^{i(k_0-j)(\pi \omega_1-\lambda)}$. Then 
\begin{eqnarray*}
\int_{0}^{s+2\pi}\tilde{h}^p_j(\tau)\text{d}\tau&=&
\int_{0}^{2\pi}\tilde{h}^p_j(\tau)\text{d}\tau+\int_{2\pi}^{s+2\pi}\tilde{h}^p_j(\tau)\text{d}\tau\\
&=&\int_{0}^{2\pi}\tilde{h}^p_j(\tau)\text{d}\tau+K^{-1}\int_{0}^{s}\tilde{h}^p_j(\tau)\text{d}\tau,
\end{eqnarray*}
and by \eqref{formulee_j} 

%$$\int_{0}^{2\pi}\tilde{h}_j\e^{i\frac{\lambda}{2\pi}(k_0-j)\tau}\text{d}\tau=0.$$
%Denote by $K=\e^{i(k_0-j)(\pi \omega_1-\lambda)}$, therefore by \eqref{formulee_j}
\begin{eqnarray}\label{ejperio}
\lefteqn{e^p_{j}(s+2\pi)=K\e^{-i\frac{1}{2}(k_0-j)(\omega_1-\frac{\lambda}{\pi})s}
\left(C^p_j-i\int_{0}^{s+2\pi}\tilde{h}^p_j(\tau)\text{d}\tau\right)
}\nonumber\\
&=&\e^{-i\frac{1}{2}(k_0-j)(\omega_1-\frac{\lambda}{\pi})s}
\left(KC^p_j-iK\int_{0}^{2\pi}\tilde{h}^p_j(\tau)\text{d}\tau-i \int_{0}^{s}\tilde{h}^p_j(\tau)\text{d}\tau\right).
\end{eqnarray}
Notice that $K\not=1$, as $\lambda\not\in \pi \Q$ and choose
$$C^p_j=\frac{iK}{K-1}\int_{0}^{2\pi}\tilde{h}^p_j(\tau)\text{d}\tau,$$
then according to \eqref{formulee_j}and \eqref{ejperio}, the function $e^p_j$ is $2\pi-$periodic.\\
Now, we show that the constants $C^p_j$ are uniformly bounded in $j\geq 0$, so that the function $v_p$ given
by \eqref{defvp} is well defined. We first need the

\begin{lem}\label{lemfourier}
Let $(h^p_j)_{j\geq 0}\in l^2(\N)$ be the family of $2\pi-$periodic functions defined by
 \eqref{defhj} and $h^p_j(s)=\sum_{n\in \Z}c^p_{l,j}\e^{ils}$ its Fourier decomposition. Then for all 
$n_1,n_2\in \N$ there exists $C^p>0$ such that for all $j\in \N$
$$\sum_{l\in \Z}j^{2n_1}l^{2n_2}|c^p_{l,j}|^2\leq C^p. $$
\end{lem}
\dem 
Consider the function $F_p\in \mathcal{C}^{\infty}\left([0,2\pi],\mathcal{S}(\R)\right) $ which defines 
the family $(g^p_j(s))_{j\geq 0}\in l^2(\N)$ with \eqref{defF}. Denote by 
 $H=-\frac12\partial^2_x+\frac12x^2$. Let $n_1, n_2 \in \N$ and decompose the function 
$\partial_s^{n_2}H^{n_1}F_p$ on the basis $(\varphi_j)_{j\geq 0}$
$$ \partial_s^{n_2}H^{n_1}F_p=\sum_{j\geq 0}\tilde{g}^p_j(s)\varphi_j(s) $$
where $(\tilde{g}^p_j)_{j\geq 0}$ is a smooth family of functions in  $l^2(\N)$.\\
Using that $H\varphi_j=(j+\frac12)\varphi_j$ and that
$F_p\in\mathcal{C}^{\infty}\left([0,2\pi],\mathcal{S}(\R)\right) $, we have for all $n_1,n_2\in \N$
%$$H^{n_1}F_p=\sum_{j\geq 0}(j+\frac12)^{n_1}{g}_j(s)\varphi_j(y) $$
%and for all $n_2\in \N$
$$\partial_s^{n_2}H^{n_1}F_p=\sum_{j\geq 0}(j+\frac12)^{n_1}(g^p_j)^{(n_2)}(s)\varphi_j(y). $$
By uniqueness of such a decomposition, 
$$\Big((j+\frac12)^{n_1}(g^p_j)^{(n_2)}\Big)_{j\geq0}=(\tilde{g}^p_j)_{j\geq0}\in l^2(\N)$$
%${g}^{(n_2)}_j(s)=\tilde{g}\in l^2(\N)$.\\
Then by the definition \eqref{defhj} of $h^p_j$, an easy induction on $n_1,n_2\in \N$ shows that 
$\big(j^{n_1}(h^p_j)^{(n_2)}\big)_{j\geq 0}\in l^2(\N)$. Write the Fourier decomposition of $h^p_j$
$$h^p_j(s)=\sum_{n\in \Z}c^p_{l,j}\e^{ils}$$ and by Parseval
$$\sum_{j\geq 0}\sum_{l\in \Z}j^{2n_1}l^{2n_2}(s)|c^p_{l,j}|^2=
\sum_{j\geq 0}j^{2n_1}\int_0^{2\pi}|(h^p_j)^{(n_2)}(s)|^2\text{d}s\leq C^p.$$
%$$j^{n_1}h_j^{(n_2)}(s)=\sum_{n\in \Z}c_{l,j}j^{n_1}(il)^{n_2}(s)\e^{ils}$$ 
%$$\sum_{j\geq 0}j^{2n_1}|h_j^{(n_2)}(s)|^2\leq C$$
%$$\sum_{j\geq 0}\sum_{l\in \Z}j^{2n_1}|l|^{2n_2}(s)|c_{l,j}|^2\leq C$$
In particular, for all $j\in \N$
$$\sum_{l\in \Z}j^{2n_1}l^{2n_2}|c^p_{l,j}|^2\leq C^p, $$
hence the result.
\findem
\\[10pt]
\noindent \textit{End of the proof of Proposition \ref{propE_p}: }
Using the Fourier decomposition of ${h}_j$ we obtain
\begin{eqnarray}\label{estCj}
C^p_j&=&\frac{iK}{K-1}\int_{0}^{2\pi}\tilde{h}^p_j(\tau)\text{d}\tau\nonumber\\
   &=&\frac{iK}{K-1}\sum_{l\in \Z}c^p_{l,j}
\int_0^{2\pi}\e^{i(l-\frac12(k_0-j)(\omega_1-\frac{\lambda}{\pi}))\tau}\nonumber\\
 &=&-i\sum_{l\in \Z}\frac{c^p_{l,j}}{l-\frac12(k_0-j)(\omega_1-\frac{\lambda}{\pi})}.
\end{eqnarray}
With Assumption \ref{assump2} we have 
\begin{eqnarray*}
\big|l-\frac12(k_0-j)(\omega_1-\frac{\lambda}{\pi})\big|&=&
\frac12|(2l-(k_0-j)\omega_1)+(k_0-j)\frac{\lambda}{\pi}|\nonumber\\
&\geq & \frac12\frac{\mu}{|(2l-(k_0-j)\omega_1,k_0-j)|^{\tau}},
\end{eqnarray*}
and for $j\geq k_0$, $|2l-(k_0-j)\omega_1|+|k_0-j|\leq 2 (|l|+|j|)$, then 
\begin{equation}\label{estCj0}
\big|l-\frac12(k_0-j)(\omega_1-\frac{\lambda}{\pi})\big|
\geq \frac{\mu'}{(|l|+|j|)^{\tau}}.
\end{equation}
Hence, from \eqref{estCj} and  \eqref{estCj0} we deduce
\begin{equation}\label{estCj1}
|C^p_j|\lesssim \sum_{l\in \Z}|c^p_{l,j}|(|j|+|l|)^{\tau}\lesssim
 \sum_{l\in \Z}|c^p_{l,j}|(|j|^{\tau}+|l|^{\tau}).
\end{equation}
\noindent By Cauchy-Schwarz and Lemma \ref{lemfourier}, from \eqref{estCj1} we obtain 
\begin{eqnarray}\label{estCj2}
|C_j|&\lesssim & \sum_{l\in \Z}\frac{1+|l|}{1+|l|}|c^p_{l,j}|(|j|^{\tau}+|l|^{\tau})\nonumber\\
      &\lesssim &\big(\sum_{l\in \Z}\frac{1}{(1+|l|)^2}\big)^{\frac12}
\big(\sum_{l\in \Z}|c^p_{l,j}|^2(1+|l|)^2(|j|^{2\tau}+|l|^{2\tau})\big)^{\frac12}\nonumber\\
&\leq & C^p.
\end{eqnarray}
Set 
$$v_p(s,x)=\sum_{j \geq 0 }e^p_j(s)w_j(s,x).$$
For all $j\in \N$, $s \longmapsto e^p_j(s)w_j(s,x)$ is continuous and there exists $c>0$ such that for 
all $j >k_0$, and for all $s\in [0,2\pi]$
$$|e^p_j(s)w_j(s,x)|\lesssim |g^p_j(s)||\varphi_j(cx)|$$
%&\lesssim& (|p_j|+|q_j|)|\varphi_j(cx)|
%\end{eqnarray*}
and this shows that $v_p \in C\left([0,2\pi],L^2(\R)\right)$. 
Now using Proposition $\ref{propregularite}$ we conclude, by uniqueness of such a solution, that 
$v_p \in \mathcal{C}^{\infty}\left([0,2\pi],\mathcal{S}(\R)\right)$.
\findem
\subsection{The nonlinear analysis and proof of Proposition \ref{propWKB}}
\begin{lem}\label{lemE_1}
The constant $E_1$ given by Proposition \ref{propE_p} writes $E_1=-\eps\delta^2C_0$ where $C_0>0$ 
 is independent of $\eps$ and $\delta$.
\end{lem}

\dem
According to formula \eqref{formuleE_p}, we only have to compute the term 
$h_{k_0}$ in the expansion \eqref{dvptQ_p}.\\
Write the expansion of $|\varphi_{k_0}|^2 \varphi_{k_0} $ on the basis $(\varphi_j)_{ j\geq0}$:
\begin{equation}\label{dvptphi_k}
|\varphi_{k_0}|^2 \varphi_{k_0}=\sum_{j\geq 0}p_j\varphi_j,
\end{equation}
with  $p_j \in \R$ and $p_j=0$ for $j-k_0=1$ $mod \;2$ as $\varphi_k(-x)=(-1)^k\varphi_k(x).$ 
Then 
\begin{eqnarray*}
|v_0|^2v_0&=&\e^{-i s E_0(k_0)}\e^{i\dot{\alpha}(s)x^2/2}\e^{-i(\frac12+k_0)\beta(s)}\e^{-\frac32\alpha(s)}|\varphi_{k_0}|^2\varphi_{k_0}\left(x\e^{-\alpha(s)}\right)\\
&=& \sum_{j\geq 0}p_j\e^{-i s E_0(k_0)} \e^{i\dot{\alpha}(s)x^2/2}\e^{-i(\frac12+k_0)\beta(s)}\e^{-\frac32\alpha(s)}    \varphi_j\left(x\e^{-\alpha(s)}\right) \\
&=& \sum_{j\geq0}f_j(s)w_j(s,x)
\end{eqnarray*}
where 
\begin{eqnarray*}
f_j(s)&=&p_j \e^{-i s( E_0(k_0)- E_0(j))} \e^{-i(k_0-j)\beta(s)}\e^{-\alpha(s)}\\
&=&p_j \e^{-i(k_0-j)(\theta(s)+\frac{s}{2}\omega_1)}\e^{-\alpha(s)}. 
\end{eqnarray*}
Therefore $f_{k_0}(s)=p_{k_0}\e^{-\alpha(s)}$ with, using \eqref{dvptphi_k},
 $p_{k_0}=\int_{\R} |\phi_{k_0}|^4>0$.\\
%which satisfies
%\begin{equation}\label{periof_j}
%f_j(s+2\pi)=\e^{-i(k_0-j)\pi \omega_1}f_j(s)=(-1)^{(j-k_0)\omega_1}f_j(s)=f_j(s)
%\end{equation}
%because $f_j=0$ for $j-k_0=1$ $mod \;2$.\\
In the same manner we write
$$x^3\varphi_{k_0}(x)=\sum_{j\geq 0}q_j\varphi_j(x)$$
with $q_j=0$ when  $j-k_0=0$ $mod \;2$ and 
\begin{eqnarray*}
R_3x^3v_0&=&R_3\e^{-i s E_0(k_0)}\e^{i\dot{\alpha}(s)x^2/2}\e^{-i(\frac12+k_0)\beta(s)}\e^{\frac52\alpha(s)}
(x\e^{-\alpha(s)})^3\varphi_{k_0}\left(x\e^{-\alpha(s)}\right)\\
&=&  \sum_{j\geq0}g_j(s)w_j(s,x)
\end{eqnarray*}
where 
\begin{eqnarray*}
g_j(s)&=&q_jR_3(s) \e^{-i s( E_0(k_0)- E_0(j))} \e^{-i(k_0-j)\beta(s)}\e^{3\alpha(s)}\\
&=&q_jR_3(s) \e^{-i(k_0-j)(\theta(s)+\frac{s}{2}\omega_1)}\e^{3\alpha(s)}. 
\end{eqnarray*}
Then $g_{k_0}=0$ and $h_{k_0}=\frac12\eps\delta^2f_{k_0}=\frac12\eps\delta^2p_{k_0}\e^{-\alpha(s)}$.\\
Finally, from \eqref{formuleE_p} we deduce 
$$E_1=-\frac{1}{4\pi}\eps\delta^2p_{k_0}\int_{0}^{2\pi}\e^{-\alpha(\tau)}\text{d}\tau=-\eps\delta^2 C_0$$
where $C_0>0$ as $p_{k_0}>0$.
\findem
\\[10pt]
\noindent Now we prove Proposition \ref{propWKB}
% Thanks to the $v_j$'s we are now able to prove Theorem $\ref{thmWKB}$, which  will be a consequence 
%of the following Proposition.
\\[10pt]
\noindent \textit{Proof of Proposition \ref{propWKB}: }
To begin with, as $\chi$ is an even function and by construction of the $v_j$'s, 
$u_p(s+2\pi,x)=u_p(s,\omega x)$.\\
Set $v(s,x)=\left(v_0+h^{\frac12}v_1+\cdots+h^{\frac{p}{2}}v_p\right)(s,x)$, then 
plugging in $(\ref{equationv})$ we obtain that the coefficient of  $h^{\frac{j}{2}}$ 
cancel for $0 \leq j\leq p$, this leads to an error term $h^{\frac{p+1}{2}}F_h(s,x)$ 
with $F_h \in \mathcal{C}^{\infty}\left([0,2\pi]\times
 ]-\frac{r_0}{\sqrt{h}},\frac{r_0}{\sqrt{h}}[,\R\right)$. Note that the function $F_h$ isn't 
defined for $x\in \R$, as it depends on $a$ which is only defined for $r\in ]-r_0,r_0[$.\\
Set $\tilde{u}_p=\delta \e^{\frac{i}{h}s}v(s,\frac{r}{\sqrt{h}})$, then coming back to 
equation $(\ref{elliptique})$, there exists  $G_h \in \mathcal{C}^{\infty}\left([0,2\pi]\times
 ]-\frac{r_0}{\sqrt{h}},\frac{r_0}{\sqrt{h}}[,\R\right)$ so that 
$$-\Delta \tilde{u}_p - \lambda_p \tilde{u}_p +\eps |\tilde{u}_p|^2\tilde{u}_p=\delta h^{\frac{p-1}{2}}h^{-\frac14}G_h. $$
(Here we loose a power of $h$ as we had multiplied by $h$ to obtain  $(\ref{equationv})$.)\\
Now 
\begin{eqnarray}
-\Delta u_p = \lambda_p u_p +\eps |u_p|^2u_p &=& 
\chi(-\Delta \tilde{u}_p - \lambda_p \tilde{u}_p)+\eps\chi^3|\tilde{u}_p|^2\tilde{u}_p \nonumber\\
&&-\chi' (2\tilde{u}_p +\frac{\partial_ra}{a}\tilde{u}_p )  -\chi'' \tilde{u}_p \nonumber\\
 &=&  h^{\frac{p-1}{2}}h^{-\frac14}\chi G_h+\eps \chi(\chi^2-1)|\tilde{u}_p|^2\tilde{u}_p  \nonumber\\
&&-\chi' (2\tilde{u}_p +\frac{\partial_ra}{a}\tilde{u}_p ) -\chi'' \tilde{u}_p \label{equationerreurG}\\
 &:=&  h^{\frac{p-1}{2}}g_p(h). \nonumber
\end{eqnarray}
First, observe that 
\begin{eqnarray}
h^{-\frac12}\int_{[0,2\pi]\times \R}\chi^2(r)\left|G_h(s,{\frac{r}{\sqrt{h}}})\right|^2\text{d}s \text{d}r
&=& \int_{[0,2\pi]\times \R}\chi^2(\sqrt{h}x)\left|G_h(s,x)\right|^2\text{d}s \text{d}x\nonumber \\
& \sim& 1.\label{estrest1}
\end{eqnarray}
Then, as $v$ is rapidly decreasing in $x$, $\tilde{u}_p$ is localized near $r=0$ in an 
interval of length $\sim \sqrt{h}$. But each of the terms $\chi(\chi^2-1)$, $\chi'$ and $\chi''$ vanish 
when $|r|\leq r_0/2$. Thus, for all $N \in \N$ there exists $C_N$ such that
\begin{eqnarray*}
\| \chi(\chi^2-1)|\tilde{u}_p|^2\tilde{u}_p\|_{L^2(M)} \leq C_N h^N, \\
\| \chi' (2\tilde{u}_p +\frac{\partial_ra}{a}\tilde{u}_p )    \|_{L^2(M)} \leq C_N h^N,\\
\| \chi'' \tilde{u}_p    \|_{L^2(M)} \leq C_N h^N.
\end{eqnarray*}
These estimates together with $(\ref{estrest1})$ yield 
$$ \| g_p(h)    \|_{L^2(M)} \lesssim \delta .  $$
Derivating a term of the right hand side of $(\ref{equationerreurG})$ costs at most $h^{-1}$ 
(when you derivate $\exp{i\frac{s}{h}}$ is variable $s$). Hence, for all $n \in \N$
\begin{equation*}
  \| g_p(h)    \|_{L^2(M)} \lesssim \delta h^{-n}  .     
\end{equation*}
\findem

\begin{lem}\label{lemreel}
Let $p\geq 1$ and  $E_p$  given by Proposition $(\ref{propE_p})$. Then 
$E_p\in \R$.
\end{lem}

\dem
 We already know that $E_0,E_1\in\R$. Let $p\geq 3$. Multiply $(\ref{exprelliptique_p})$ by $\overline{u}_p$, integrate on $M$ and take the imaginary part
$$0=\text{Im} \lambda_p \| u_p \|^2_{L^2}+ h^{\frac{p-1}{2}} \text{Im} \int g_p(h)\overline{u}_p. $$
As $\| u_p \|_{L^2}\sim 1$ and $\| g_p \|_{L^2}\lesssim 1$, we obtain the estimate
$$|\text{Im}\lambda_p|\lesssim  h^{\frac{p-1}{2}}\| g_p \|_{L^2}\| u_p \|_{L^2} \lesssim  h^{\frac{p-1}{2}}$$
and as 
$$\text{Im}\lambda_p=-2(\text{Im} E_2+h^{\frac12}\text{Im} E_2+\cdots+ h^{\frac{p-1}{2}}\text{Im} E_p)  $$
it follows that for all $0\leq j\leq p-1$, $\text{Im} E_j=0$, i.e. $E_j\in\R$.
\findem
%\end{document}

\section{The instability for the nonlinear Schr\"odinger equation}
\subsection{The error estimate}
\begin{prop}\label{properreur}
Let $\alpha>0$ and let $v \in H^2(M)$ be such that 
$$ \| v\|_{L^2} \lesssim 1,\;\;\| v\|_{L^{\infty}} \lesssim h^{-\frac{1}{4}+\sigma},\;\;
\|\Delta  v\|_{L^{\infty}} \lesssim h^{-\frac{9}{4}+\sigma}, $$ 
and suppose that $v$ satisfies
\begin{equation*}
i \partial_t v+\Delta v   =  \eps|v|^{2} v+ h^{\alpha}R(h)
\end{equation*}
with for all $\beta \in [0,2]$,  $\|R(h)\|_{H^{\beta}} \lesssim h^{-\beta}$.
Let $u$ be solution of
\begin{equation*}
\left\{
\begin{array}{l}
i \partial_t u+\Delta u=\eps|u|^{2} u\\
u(0,x)=v(0,x).
\end{array}
\right.
\end{equation*}
Then, if $\alpha>\frac14+3\sigma$ we have
$$\|(u-v)(t_h)\|_{H^{\sigma}}\longrightarrow0\quad \text{when}\quad h\longrightarrow0,$$
where $t_h\sim h^{\frac{1}{2}-2\sigma}\log(\frac1h)$.
\end{prop}

\dem
Define $w=u-v$ and 
\begin{equation*}
E(t)=\|w\|^2_{L^2}+\|h^2 \Delta w\|^2_{L^2}.
\end{equation*}
We have $E(0)=0$ and the following estimates:
\begin{equation}\label{estenergie}
\|w\|_{L^2}  \leq  E^{\frac{1}{2}}, \quad \|\Delta w\|_{L^2}  \leq  h^{-2}E^{\frac{1}{2}},
\quad \|\nabla w\|_{L^2}  \leq  h^{-1}E^{\frac{1}{2}}.
\end{equation}
The function $w$ satisfies the equation
\begin{equation}\label{err}
i \partial_t w+\Delta w   =  \eps(|w+v|^{2} (w+v)-|v|^2v)+R(h).
\end{equation}
The energy method gives
\begin{eqnarray*}
\frac12\frac{\text{d}}{\text{d}t}\|w\|^2_{L^2}& =& \text{Im} \int\overline{w}
\left(\eps(|w+v|^{2} (w+v)-|v|^2v)+R(h)\right) \nonumber \\
 & \lesssim & h^{\alpha}\|w\|_{L^2}+\|w\|^4_{L^4}+\|w\|^2_{L^2}\|v\|^2_{L^{\infty}}.
\end{eqnarray*}
The Gagliardo-Nirenberg inequality gives
$$\|w\|^4_{L^4} \lesssim \|w\|^2_{L^2}\|\nabla w\|^2_{L^2} \lesssim h^{-2}E^2,$$
and as $\|v\|_{L^{\infty}}\lesssim h^{-\frac{1}{4}+\sigma}$, we obtain
\begin{equation}\label{estwL2}
\frac{\text{d}}{\text{d}t}\|w\|^2_{L^2} \lesssim  h^{\alpha}E^{\frac{1}{2}}+ 
h^{-\frac{1}{2}+2\sigma}E+h^{-2}E^2.
\end{equation}
Now, apply $\Delta$ to $(\ref{err})$
\begin{equation}\label{deltaerr}
i \partial_t \Delta w+\Delta^2 w   = \eps \Delta A+\Delta R(h)
\end{equation}
with
\begin{eqnarray*}
A &=& |w+v|^{2} (w+v)-|v|^2v\\
 &=& 2w|v|^2+\overline{w}v^2+w^2\overline{v}+2|w|^2v+|w|^2w,
\end{eqnarray*}
then
\begin{eqnarray*}
|\Delta(A)| &\lesssim & |v|^2|\Delta w|+|v||\nabla v||\nabla w|+|\nabla v|^2|w|+|v||\Delta v||w| \nonumber \\
     & & +|\Delta v||w|^2+|w|^2|\Delta w|+|w||\nabla w|^2
\end{eqnarray*}
hence
\begin{eqnarray}
\|\Delta(A) \|_{L^2} &\lesssim & \|v\|_{L^{\infty}}^2\|\Delta w\|_{L^2}+\|v\|_{L^{\infty}}\|\nabla v\|_{L^{\infty}}\|\nabla w\|_{L^2}+\|\nabla v\|_{L^{\infty}}^2\|w\|_{L^2}\nonumber \\
& & +\|v\|_{L^{\infty}}\|\Delta v\|_{L^{\infty}}\|w\|_{L^2}
  +\|\Delta v\|_{L^{\infty}}\|w\|_{L^4}^2 \label{estA0}\\
& &+\|w\|_{L^{\infty}}^2\|\Delta w\|_{L^2}+\|w\|_{L^2}\|\nabla w\|_{L^4}^2.\nonumber
\end{eqnarray}
The following inequality holds in dimension $2$
\begin{equation*}
\|w\|_{L^{\infty}}\lesssim \|w\|_{L^2}^{\frac12}\|\Delta w\|_{L^2}^{\frac12} \lesssim h^{-1} E^{\frac12},
\end{equation*}
and with $(\ref{estenergie})$ and $(\ref{estA0})$ we deduce
\begin{equation*}
\|\Delta(A) \|_{L^2}\lesssim h^{-\frac52+2\sigma}E^{\frac12}+ h^{-\frac{13}{4}+\sigma}E 
+ h^{-4} E^{\frac32}.
\end{equation*}
But $$h^{-\frac{13}{4}+\sigma}E=h^{-\frac{5}{4}+\sigma}E^{\frac14} h^{-2} E^{\frac34}\lesssim
 h^{-\frac52+2\sigma}E^{\frac12}+  h^{-4} E^{\frac32},$$
and we obtain
\begin{equation}\label{estA}
\|\Delta(A) \|_{L^2}\lesssim h^{-\frac52+2\sigma}E^{\frac12}+ h^{-4} E^{\frac32}.
\end{equation}
Now, using $(\ref{estA})$ and $\|\Delta(R(h)) \|_{L^2}\lesssim h^{\alpha-2}$, the energy method and 
the Cauchy-Schwarz inequality gives
\begin{eqnarray}
\frac12\frac{\text{d}}{\text{d}t}\|\Delta w\|^2_{L^2}& =& \text{Im} \int\Delta\overline{w}
\left(A+R(h)\right) \nonumber \\
 & \lesssim & h^{-2} E^{\frac12}(h^{\alpha-2}+ h^{-\frac52+2\sigma}E^{\frac12}+ h^{-4} E^{\frac32}),\label{estdeltawL2}
\end{eqnarray}
therefore from $(\ref{estwL2})$ and $(\ref{estdeltawL2})$ we have
\begin{equation*}
\frac{\text{d}}{\text{d}t} E \lesssim  h^{\alpha}E^{\frac{1}{2}}+ 
h^{-\frac{1}{2}+2\sigma}E+h^{-2}E^2.
\end{equation*}
Interpolation gives
\begin{equation*}
\|w\|_{H^{\sigma}} \lesssim \|w\|_{L^2}+\|w\|_{\dot{H}^{\sigma}} 
\lesssim \|w\|_{L^2} +  \|w\|^{1-\frac{\sigma}{2}}_{L^2} \|\Delta w\|^{\frac{\sigma}{2}}_{L^2}
\lesssim h^{-\sigma}E^{\frac12}:=F
\end{equation*}
and $F$ satisfies $F(0)=0$ and

\begin{equation}\label{eqF}
\frac{\text{d}}{\text{d}t} F \lesssim  h^{-\sigma+\alpha} +h^{-\frac{1}{2}+2\sigma}F 
+h^{-2+2\sigma}F^3.
\end{equation}
As long as $h^{-2+2\sigma}F^3 \lesssim h^{-\frac{1}{2}+2\sigma}F $, i.e. $F \lesssim h^{\frac34}$, 
we can write 
\begin{equation*}
\frac{\text{d}}{\text{d}t} F \lesssim  h^{-\sigma+\alpha} +h^{-\frac{1}{2}+2\sigma}F
\end{equation*}
and the Gronwall inequality yields 
\begin{equation*}
F \lesssim h^{\alpha+\frac12-3\sigma} \e^{C+h^{-\frac{1}{2}+2\sigma}t}.
\end{equation*}
The non linear term in $(\ref{eqF})$ can be removed with the continuity argument for times such that
\begin{equation*}
h^{\alpha+\frac12-3\sigma }\e^{C+h^{-\frac{1}{2}+2\sigma}t}\lesssim h^{\frac34+\eta}
\end{equation*}
with $\eta>0$
i.e. for $t \lesssim (\alpha-\frac14-3\sigma-\eta)h^{-\frac{1}{2}+2\sigma}\log{\frac{1}{h}}$, which is 
possible with $\eta$ small enough as we assume $\alpha>\frac14+3\sigma$.
\findem

\begin{cor}\label{coransatz} Let  $\kappa>0$, $0<\sigma<\frac14$ and set $\delta =\kappa h^{\sigma}$. Denote by $v=\e^{-i\lambda_3 t}u_3$ 
where $u_3$ and $\lambda_3$ are defined by  $(\ref{expru_p})$ and $(\ref{exprlambda_p})$ respectively.\\
Let $u$ be solution of
\begin{equation*}
\left\{
\begin{array}{l}
i \partial_t u+\Delta u=\eps|u|^{2} u\\
u(0,x)=v(0,x).
\end{array}
\right.
\end{equation*}
Then $\|v\|_{H^{\sigma}}\sim 1$ and 
$$\|(u-v)(t_h)\|_{H^{\sigma}}\longrightarrow0\quad \text{when}\quad h\longrightarrow0,$$
where $t_h\sim h^{\frac{1}{2}-2\sigma}\log(\frac1h)$. 
\end{cor}

\dem
The result directly follows from Propositions  $(\ref{propWKB})$ and  $(\ref{properreur})$, as 
for all $0<\sigma<\frac14$, we have $\sigma+1 >\frac14+3\sigma$.
\findem

\subsection{The instability argument}
Let $\kappa,\kappa_h>0$ and consider $v=v^1$ defined in Corollary $\ref{coransatz}$ associated with $\kappa$ 
and $v^2$ associated with $\kappa_h$. Let $u$ be a solution of 
\begin{equation*}
\left\{
\begin{array}{l}
i \partial_t u^j+\Delta u^j=\eps|u^j|^{2} u^j\\
u^j(0,x)=v^j(0,x)
\end{array}
\right.
\end{equation*}
and $t_h\sim h^{\frac{1}{2}-2\sigma}\log(\frac1h)$. Then
\begin{eqnarray}
\|(u^2-u^1)(t_h)\|_{H^{\sigma}}&\geq & \|(v^2-v^1)(t_h)\|_{H^{\sigma}}-\|(u^2-v^2)(t_h)\|_{H^{\sigma}}\nonumber \\
&&-\|(u^1-v^1)(t_h)\|_{H^{\sigma}}.\label{comparaison}
\end{eqnarray}
From Corollary $\ref{coransatz}$ we deduce that for $j=1,2$
\begin{equation}\label{limite}
 \|(u^j-v^j)(t_h)\|_{H^{\sigma}}\longrightarrow0.
\end{equation}
Observe that
$$\|(v^2-v^1)(t_h)\|_{H^{\sigma}}\sim \left|\e^{-i\lambda_3^2 t_h} - \e^{-i\lambda_3^1 t_h}\right|
=\left| \e^{i(\lambda_3^2- \lambda_3^1)t_h}-1\right|,$$
from Proposition  $\ref{lemE_1}$ we have 
$$(\lambda_3^2- \lambda_3^1)t_h\sim h^{2\sigma -1}(\kappa-\kappa_h)t_h\sim (\kappa-\kappa_h)\log{\frac{1}{h}}.$$
It is possible to choose $\kappa_h$ such that $\kappa_h \longrightarrow \kappa$ and 
$(\kappa-\kappa_h)\log{\frac{1}{h}}\longrightarrow \infty$. Then using $(\ref{comparaison})$ and $(\ref{limite})$
$$\limsup_{h\longrightarrow 0}\|(u^2-u^1)(t_h)\|_{H^{\sigma}}\geq\limsup_{h\longrightarrow 0}
 \|(v^2-v^1)(t_h)\|_{H^{\sigma}}\geq2 $$
even though
$$\|(u^2-u^1)(0)\|_{H^{\sigma}}= \|(v^2-v^1)(0)\|_{H^{\sigma}}\sim |\kappa-\kappa_h|$$
which tends to $0$ with $h$. According to Definition $\ref{definstab}$, we have proved Proposition $\ref{corinstab}$.

\nocite{Helffer}
\nocite{CCT1}
\nocite{CCT2}

%Pour compiler:
%latex doc.tex
%latex doc.tex
%bibtex doc
%latex doc.tex
%

%\begin{thebibliography}{10}

\bibliographystyle{plain}
%\bibliography{biblio}
%\email{laurent.thomann@math.u-psud.fr}

\vspace{1cm}
\noindent
{\sc{L. Thomann, Universit\'e Paris-Sud, Math\'ematiques, B\^at 425, 91405
Orsay Cedex.}}\\
{\it{E-mail}}:\textrm{ laurent.thomann@math.u-psud.fr}

\end{document}